\documentclass[11pt,a4paper]{amsart}
\usepackage{amssymb}
\usepackage{amsmath}
\usepackage{amsfonts}

\newcommand{\Q}{\mathbb{Q}}

\newcommand{\F}{\mathbb{F}}

\newtheorem{dummy}{Dummy}

\newtheorem{lemma}[dummy]{Lemma}
\newtheorem{theorem}[dummy]{Theorem}

\theoremstyle{definition}

\theoremstyle{remark}

\begin{document}
\bibliographystyle{amsalpha}
%%%%%%%%%%%%%%%%%%%%%%%%%%%%%%%%%%%%%%%%%%%%%%%%%%%%%%%%%%%%%%%%%%%%%%%
%%%%%%%%%%%%%%%%%%%%%%%%%%%%%%%% Title %%%%%%%%%%%%%%%%%%%%%%%%%%%%%%%%
\author{Sandro Mattarei}
\email{mattarei@science.unitn.it}
\urladdr{http://www-math.science.unitn.it/\~{ }mattarei/}
\address{Dipartimento di Matematica\\
  Universit\`a degli Studi di Trento\\
  via Sommarive 14\\
  I-38050 Povo (Trento)\\
  Italy}

\title{Inverse-closed additive subgroups of fields}

\begin{abstract}
We describe the additive subgroups of fields which are closed with respect to taking inverses.
In particular, in characteristic different from two any such subgroup is either a subfield
or the kernel of the trace map of a quadratic subextension of the field.
\end{abstract}

%\date{18 November 2005}



\thanks{Partially supported by MIUR-Italy via PRIN 2003018059
 ``Graded Lie algebras and pro-$p$-groups: representations,
 periodicity and derivations''.}

\maketitle

\thispagestyle{empty}
%%%%%%%%%%%%%%%%%%%%%%%%%%%%%%%%%%%%%%%%%%%%%%%%%%%%%%%%%%%%%%%%%%%%%%%
%%%%%%%%%%%%%%%%%%%%%%%%%%%%%%%%%%%%%%%%%%%%%%%%%%%%%%%%%%%%%%%%%%%%%%%

%----------------------------------------------------------------
\section{Introduction}\label{sec:intro}

The following result of Hua,
as stated in~\cite[Theorem~1.15]{Artin:geometric},
plays a role in connection with the fundamental theorem of projective geometry
(see~\cite[Chapter~II, Sections~9 and~10]{Artin:geometric}):
an additive map between division rings sending $1$ to $1$ and inverse elements to inverse elements
is either an isomorphism or an antiisomorphism of rings.
The first step in the proof is showing that the map preserves the operation $(a,b)\mapsto aba$.
This follows from Hua's identity
(first mentioned in~\cite{Hua:sfield_properties}, but see~\cite[page 2]{Jac:Jordan} or~\cite[page 89]{Jac:BAI}
for the more manageable form given here)
\begin{equation}\label{eq:Hua}
a-(a^{-1}+(b^{-1}-a)^{-1})^{-1}=aba,
\end{equation}
which holds in any associative ring provided all inverses involved are defined, that is,
provided $a,b$ and $ab-1$ are invertible.
The rest of the proof (originally given in~\cite{Hua:sfield_automorphisms})
does not use that the map preserves inverses, but rather the product $aba$.
This part of the proof has been later generalized in a number of directions,
notably to arbitrary domains by Jacobson and Rickart,
see the references given in~\cite[page~3]{Jac:Jordan}.

A problem of a similar flavour as Hua's result,
but which seems not to have received attention,
is a description of the additive subgroups
of division rings which contain the inverses of their nonzero elements.
In the present note we fill this gap in the commutative case.

For any subset $S$ of a field $E$ we write $S^{-1}=\{s^{-1}\mid 0\neq s\in E\}$.
We call $S$ {\em inverse-closed} if $S^{-1}\subseteq S$.
We prove the following results.

\begin{theorem}\label{thm:odd}
Let $E$ be a field of characteristic different from two
and let $A$ be a non-trivial inverse-closed additive subgroup of $E$.
Then $A$ is either a subfield of $E$ or the set of elements of trace zero in some quadratic field extension
contained in $E$.
\end{theorem}

Conversely, it is plain that the set of elements of trace zero in any quadratic field extension
contained in $E$ is inverse-closed.

\begin{theorem}\label{thm:even}
Let $E$ be a field of characteristic two
and let $A$ be an inverse-closed additive subgroup of $E$.
Then $A$ is an $F^2$-subspace of $F$ for some subfield $F$ of $E$.
\end{theorem}

Conversely, any $F^2$-subspace of a subfield $F$ of characteristic two is clearly inverse-closed.

Because of the method of proof employed, based in particular on Hua's identity, it is natural to
ask for extensions of Theorems~\ref{thm:odd} and~\ref{thm:even} to division algebras.
However, the correct extension is not yet clear to this author from a variety of examples found.

The inversion map in finite fields is of cryptographic interest.
For example, inversion in the finite field of $2^8$ elements
is the nonlinear transformation employed in the S-boxes
in the Advanced Encryption Standard (Rijndael)~\cite{AES}.
In view of possible applications, as in~\cite{CDVSV}, the special case of Theorems~\ref{thm:odd} and~\ref{thm:even}
where $E$ is a finite field deserves the following separate mention.

\begin{theorem}\label{thm:finite}
Let $E$ be a finite field
and let $A$ be a non-trivial inverse-closed additive subgroup of $E$.
Then $A$ is either a subfield of $E$ or the set of elements of trace zero in some quadratic field extension
contained in $E$.
\end{theorem}

Note that when $E$ has characteristic two the two alternatives in the conclusion coincide.
The finiteness assumption on $E$ allows various proofs of Theorem~\ref{thm:finite}
which differ from that of the general case.
For example, one of the advantages of finite fields is that an arbitrary subset can be described by the unique
monic polynomial which has the elements of the subset as simple roots.
Algebraic properties of the subset often translate into properties of the corresponding polynomial.
%We record a particularly short one based on $p$-polynomials.
In particular, here we record a proof of Theorem~\ref{thm:finite} based on $p$-polynomials.

I am grateful to A.~Caranti for asking the question answered in Theorem~\ref{thm:finite}.

%----------------------------------------------------------------
\section{Proofs}

Note that an inverse-closed additive subgroup $A$ of $E$ is necessarily
a subspace of $E$ over its prime field.
This is clear if $E$ has positive characteristic.
If $E$ has characteristic zero, then $A$ is a $\Q$-subspace of $E$, because $(mn^{-1})a=m(na^{-1})^{-1}\in A$
for $a\in A^\ast$ and $m,n$ integers with $n\neq 0$.

\begin{lemma}\label{lemma:Hua}
Let $A$ be an inverse-closed additive subgroup of $E$.
Then $a^2b\in A$ for all $a,b\in A$.
Furthermore, if $E$ has characteristic different from two, then $abc\in A$ for all $a,b,c\in A$.
\end{lemma}

\begin{proof}
Hua's identity~\eqref{eq:Hua} implies that
$aba=a^2b\in A$
for all $a,b\in A$,
the degenerate cases where one or more of $a,b$ and $ab-1$ vanish being obvious.
The second assertion follows from the identity
$2abc=(a+c)^2b-a^2b-c^2b$.
\end{proof}

It follows at once by taking $c=1$ in Lemma~\ref{lemma:Hua}, that the only inverse-closed additive subgroups containing $1$
of a field $E$ of characteristic not two are the subfields of $E$.
According to Theorem~\ref{thm:even}, this assertion
does not extend to arbitrary fields of characteristic two.

\begin{proof}[Proof of Theorem~\ref{thm:odd}]
The inverse-closed subset
$K=\{ab\mid a,b\in A\}$
of $E$ is a subring, and hence a subfield, because
$ab-cd=a(b-a^{-1}cd)\in K$ and
$(ab)(cd)=(abc)d\in K$
for all $a,b,c,d\in A$ with $a\neq 0$.
Choose a nonzero element $a\in A$.
Then $Aa^{-1}\subseteq K$, and $Fa\subseteq A$
because of the second assertion of Lemma~\ref{lemma:Hua}.
Hence $A=Fa$, with $a^2\in K$.
We conclude that either $A$ coincides with the subfield $K$ of $E$, or is the set of
elements of trace zero in the quadratic field extension $K(a)=K+Ka$ of $E$.
\end{proof}

The following example shows that the subset $\{ab\mid a,b\in A\}$ of $E$,
which we have used in the proof of Theorem~\ref{thm:odd},
need not be a subfield when $E$ has characteristic two.
Let $E=\F_2(u_1,u_2,u_3,u_4)$ be a purely transcendental extension of transcendence degree four of the field
of two elements $\F_2$, and let $A=E^2u_1+E^2u_2+E^2u_3+E^2u_4$.
Then
$\{ab\mid a,b\in A\}=E^2+\sum_{i<j}E^2u_iu_j$
is a vector space over $E^2$ of dimension $7$, and hence not a subfield of $E$.

\begin{proof}[Proof of Theorem~\ref{thm:even}]
Let $R$ be the subring generated by the squares of the elements of $A$,
and let $K$ be the subfield of $E$ generated by $R$.
The first assertion of Lemma~\ref{lemma:Hua} shows inductively that
$A$ is an $R$-submodule of $E$.
Since $ar^{-1}=(a^{-1}r)^{-1}\in A$ for $a\in A$ and $r\in R$
we conclude that $A$ is a $K$-subspace of $E$.
Finally, if $F$ is the subfield of $E$ generated by $A$ then $K=F^2$ and $A\subseteq F$, as desired.
\end{proof}

We conclude by giving a proof of Theorem~\ref{thm:finite} based on $p$-polynomials.
A $p$-polynomial, over a field of positive characteristic $p$,
is a polynomial all whose monomials have exponents equal to powers of $p$.
A basic property of $p$-polynomials is that their sets of roots in any field are additive subgroups.
We refer to Chapter~3 of~\cite{LN} for an extensive discussion of $p$-polynomials.
We also need the concept of self-reciprocal polynomial.
For a polynomial $f(x)=\sum_{i=0}^na_ix^i$ with $a_0a_n\neq 0$ we define its {\em reciprocal} polynomial as
$x^nf(1/x)=\sum_{i=0}^na_{n-i}x^i$.
The roots of the reciprocal polynomial are clearly the inverses of the roots of the original polynomial,
with corresponding multiplicities.
We call {\em self-reciprocal} a polynomial $f(x)=\sum_{i=0}^na_ix^i$, with $a_0a_n\neq 0$,
which equals its reciprocal polynomial up to a scalar factor.
This clearly implies that $a_n=\pm a_0$.
For a polynomial with nonzero constant term and with distinct roots,
being self-reciprocal is equivalent to its set of roots being inverse-closed.

\begin{proof}[Proof of Theorem~\ref{thm:finite}]
Let $E$ have order $p^f$.
Then $E$ is the splitting field over $\F_p$ of the polynomial $x^{p^f}-x$.
According to Theorems~3.50 and~3.52 of~\cite{LN}, the additive subgroups $A$ of $E$, that is, its $\F_p$-subspaces, are
in a bijection with the monic divisors $f_A(x)$ of $x^{p^f}-x$ which are $p$-polynomials, given by letting such a polynomial
correspond to the set $A$ of its roots.
An additive subgroup $A$ of $E$ is inverse-closed if and only if $f_A(x)/x$ is self-reciprocal.
Since $f_A(x)$ is a $p$-polynomial, degree reasons easily imply that it is a binomial,
and hence has the form $x^{p^r}-x$ or $x^{p^r}+x$, for some $r$.
In the former case $A$ is the subfield of $E$ of order $p^r$.
In the latter case we may assume that $E$ has characteristic different from two, and hence $A$ is not a subfield.
Then the roots in $E$ of the polynomial $x^{p^{2r}}-x$, which include $1$, form a subfield of $E$, of order
a divisor of $p^{2r}$ larger than $|A|=p^r$.
Hence the roots of the polynomial form a subfield of $E$ of order $p^{2r}$.
The roots of $x^{p^r}+x$ form the kernel of the trace map of this subfield over its subfield of order $p^r$.
\end{proof}

\bibliography{References}

\end{document}